\documentclass[a4paper,12pt]{article}
\usepackage{amsmath,amssymb}

\begin{document}

\begin{center}     {\large A Note on the Dual of $L_1$ for general measure spaces} \end{center}\vspace{1ex}
\begin{center} Alok Goswami and B. V. Rao\end{center}
\noindent
[{\it This note is not meant for publication; it is to appear in the form of a series of exercises in a forthcoming book 
jointly authored by us.}]\\

If $(\Omega, {\cal A}, \mu)$ is a $\sigma$-finite measure space, then it is 
well-known that for every continuous linear functional 
 $T$ on $L_1$ there is unique $g\in L_\infty$  such that (i) 
$Tf = \int\!gf$ for all $f\in L_1$ and (ii) $\|T\| = \|g\|_\infty$.\\

In the book [1] a Remark on page 290 (after the 
proof of Theorem 5 in Section IV.8.5) says that this holds for more general measure spaces  
$(\Omega, {\cal A}, \mu)$, provided $\Omega$ is a disjoint union of a family 
$\{\Omega_\alpha\}$ of $\sigma$-finite sets   
 satisfying the conditions that, any set $A\in{\cal A}$ of finite $\mu$-measure intersects at most countably  
many of the sets $\Omega_\alpha$ and any $A\in{\cal A}$ with $\mu(A\cap \Omega_\alpha)=0$ 
for all $\alpha$ must have $\mu(A) =0$. By a ``$\sigma$-finite set" here is meant a set 
in ${\cal A}$, which is a countable union of sets of finite $\mu$-measure. 

The purpose of this note to take the above remark forward and arrive at a simple description 
of $L_1^{^*}$ for any general measure space. We also want to point out here that while we consider  
here only real function spaces, the same would apply to complex function spaces as well.\\

Let $(\Omega, {\cal A}, \mu)$  be any measure space and let ${\cal S}$ denote the collection of 
all $\sigma$-finite sets in ${\cal A}$. A collection $G= \{g_A: A\in{\cal S}\}$  
is going to be called a {\bf germ} if: \\
(i) $g_A$, for each $A\in{\cal S}$, is an extended real-valued 
function defined on $A$ and measurable with respect to the restricted $\sigma$-field, \\
(ii) $\|G\|=\sup \{\|g_A\|_\infty : A\in{\cal S}\}<\infty$, that is, the 
$L_\infty$ norms of all the $g_A$ are bounded and their sup is $\|G\|$, and \\ 
(iii) $g_A = g_B$ a. e. on $A\cap B$, that is the family $\{g_A : A\in{\cal S}\}$ is `consistent'.\\

In the sequel, a germ  $\{g_A: A\in{\cal S}\}$ will be denoted simply by $(g_A)$, with the understanding 
that the sets $A$ always run over ${\cal S}$ only. 
Just as we identify functions a.e., we identify germs also: $G=(g_A)$
and $H=(h_A)$ are identified, and we write  
$G\sim H$,  if $g_A=h_A \;\;a.e.$ for all $A$. 

Let us denote the space of (equivalence classes of) germs as $L_{\infty\infty}$. The following are then easy to prove. \vspace{1ex} \\
(i)  $G+H = (g_A+h_A: A\in {\cal S})$. This is well defined, that is,
if $G\sim G^*$ and $H\sim H^*$ then $G+H\sim G^*+H^*$. \vspace{.5ex} \\
(ii) $cG = (cg_A: A\in {\cal S})$. This is well defined: $G\sim G^*$ then $cG\sim cG^*$. \vspace{.5ex} \\
(iii) $L_{\infty\infty}$ is a linear space with zero element $0= (g_A)$ where \vspace{.5ex} each $g_A=0$ a.e. \\
(iv) $\|G\|$ is a norm on $L_{\infty\infty}$. \vspace{.5ex}  \\
(v) The space  $L_{\infty\infty}$ is complete.  \vspace{.5ex}  \\
 (vi) Let  $G=(g_A)$.  If $f\in L_1(\Omega, {\cal A}, \mu)$ and  $F=(f\neq 0)$;
 put  $Tf = \int fg_F$.  Then $T$ is a well defined linear functional on $L_1$ and
$\|T\|=\|G\|$. \vspace{.5ex}  \\
(vii)  Conversely, let $T$  be a given continuous linear functional on
 $L_1(\Omega, {\cal A},\mu)$. For each $A\in{\cal S}$
 applying the sigma-finite case, get $g_A$ defined on $A$ so that
 $Tf=\int fg_A$ for all
 $f\in L_1$ with support contained in $A$. 
 Let $G=(g_A: A\in {\cal S})$. Then $G$ is a germ and $\|G\|=\|T\|$. \vspace{.5ex}  \\
 (viii)  Dual of $L_1$ is the space $L_{\infty\infty}$ with its norm $\|G\|$.\\

A very natural and important question to be asked now is whether every germ is given by a function. In other 
words, given a germ $G = (g_A)$, \vspace{.5ex} \\
Is there a measurable  function $g$  such that $g=g_A$ a.e.\;on $A$ for each $A\in {\cal S}$? The answer is: in general, NO! 
The following example illustrates this. \\

Consider $\Omega = [0,1]\times [0,1]$ with ${\cal A}$ as its Borel $\sigma$-field. Let $\mu$ be the measure on ${\cal A}$ 
defined as $\mu(A) = \sum\limits_x \lambda(A_x^v) +\sum\limits_y \lambda(A_y^h)$, where $A_x^v=\{y : (x,y)\in A\}$ and 
$A_y^h=\{x : (x,y)\in A\}$ denote respectively the vertical section of $A$ at $x$ and horizontal section of $A$ at $y$. 
Of course, $\lambda$ denotes Lebesgue measure on ${\mathbb R}$. Thus, $\mu$-measure of a set $A$ is the sum of 
Lebesgue measures of all its horizontal sections and all its vertical sections.

Note that a set $A\in{\cal A}$ is $\sigma$-finite only if  $\lambda(A_x^v) =0$, for all but countably many
 $x$,  and also, $\lambda(A_y^h)=0$, for all but countably many $y$. Now, for $f\in L_1$, define 
$${\textstyle Tf = \sum_x \int\! f(x,y) y dy\;+\;  \sum_y\int\! f(x,y) xdx}$$ 
In view of the fact that $\{f\ne 0\}$ is a $\sigma$-finite set for $f\in L_1$ and the observation made above, it is clear that 
the two sums in the definition of $Tf$ are both countable sums and so make sense. Next,  $f\in L_1$ also implies that both 
the sums actually converge and indeed $|Tf|\leq \|f\|_1$. In other words, $T$ defines a bounded linear functional on $L_1$.

 We now describe the germ that $T$ is given by. As noted above, given any $\sigma$-finite set $A$, there exist 
sequences  $(x_i: i\geq 1)$ and $(y_j: j\geq 1)$ such that
 $\lambda(A_x^v)=0$ for $x\not\in(x_i: i\geq 1)$ and $\lambda(A_y^h) =0$ for 
 $y\not\in (y_j: j\geq 1)$. Define a function $g_A$ as follows: \newpage
\noindent
 $g_A(x,y) =0$, if $x\neq x_i$ for all $i$ or if $y\neq y_j$ for all $j$ or if $(x,y) = (x_i,y_j) $ for some $(i,j)$.
 For all other points, put $g_A(x,y) = y$ if $x=x_i$ for some $i$ and  $g_A(x,y) = x$ if $y=y_j$ for some $j$.

Clearly, $G=(g_A)$ defines a germ and it is easy to see that  $Tf$ is given by the germ $G$, in the sense 
described in (vi) above.
 
Is there one function $g$ on $\Omega$ that gives this germ? That is, is there a measurable function
$g$ on $\Omega$ such that $g=g_A$ a.e. for all $A\in {\cal S}$? One can verify quite easily that the answer to 
the question is: No!\\

Finally, let us consider the set up given in the remark in the book by Dunford-Schwarz as mentioned at the start.  Suppose we 
have a 
measure space $(\Omega, {\cal A}, \mu)$, where $\Omega$ is disjoint union of $\sigma$-finite sets $\{\Omega_\alpha\}$  
 satisfying just one of the conditions in that remark, namely,  \vspace{-1ex} 
$$ A\in {\cal A} \mbox { and } \mu(A\cap \Omega_\alpha)=0 \mbox { for all } \alpha \mbox { imply } \mu(A) =0 \;
\cdots  (*) \vspace{-1ex} $$  
 Let us assume, in addition, that the family $\{\Omega_\alpha\}$ also \vspace{-1ex}  satisfies the condition   
 $$\textstyle{ A\subset \Omega \mbox { and } A\cap \Omega_\alpha \in {\cal A}, \mbox { for all } \alpha, 
\mbox { imply } A\in {\cal A}\;\cdots (**)} \vspace{-1ex} $$ 
Then every germ is indeed given by a  measurable function. Given any germ $G=(g_A)$, one can just put $g=g_{\Omega_\alpha}$ 
on $\Omega_\alpha$. Condition $(*)$ implies that this is well defined. Also, if one takes equivalent germs, one gets equivalent functions.  
Finally, condition $(**)$ implies that the function $g$ defined this way on $\Omega$ is measurable.\\
 
However, if only $(*)$ holds, but $(**)$ does not hold, then it is not always possible to represent a germ by a single measurable 
function. Here is an example. On $\Omega=[0,1]\times[0,1]$, consider the family $\{\Omega_y=L_y\cap\Omega,\, y\in [0,1]\}$ 
 of subsets, where $L_y$ denotes the horizontal line passing through $(0,y)$. Let ${\cal A}$ be the $\sigma$-field on $\Omega$, 
consisting of all those Borel sets $A$ such that either $A$ or $A^c$ is contained in countably many of the $\Omega_y$. Defining 
$\mu(A)=\sum_y\lambda(A_y^h)$,  for $A\in{\cal A}$, where $A_y^h$ denotes, as before, the horizontal section of $A$ at $y$, 
one can easily see that both the conditions in the remark in Dunford-Schwarz 
hold with $\{\Omega_y\}$ as the disjoint family of $\sigma$-finite sets (in fact, each $\Omega_y$ has $\mu$-measure $1$). 
However, \vspace{-.7ex}  
$$\textstyle{Tf=\sum_y\int\! f(x,y)xdx,\;\; f\in L_1}\vspace{-.7ex} $$ 
defines a bounded linear functional on $L_1$, but the germ that $T$ is given by, cannot be represented by a single measurable function. \\

Reference:\\

[1] N. Dunford and J. T. Schwartz [1957]: Linear Operators, Part I: General Theory;
Interscience Publishers, New York
\vspace{3ex}

\noindent
Alok Goswami $<alok.gosw@gmail.com>$\\
Indian Association for the Cultivation of Science, Kolkata \vspace{1ex} \\
\noindent
B. V. Rao $<bvrao@cmi.ac.in>$\\
Chennai Mathematical Institute, Chennai.

\end{document}